\documentclass{amsart}
\usepackage{lastpage}
\usepackage{yhmath,verbatim}
\usepackage[all]{xy}

\newcommand{\bdis}{\begin{displaymath}}
\newcommand{\edis}{\end{displaymath}}
\newcommand{\be}{\begin{equation}}
\newcommand{\ee}{\end{equation}}
\newcommand{\mbb}{\mathbb}
\newcommand{\mcal}{\mathcal}

\newcommand{\vp}{\varphi}

\newcommand{\zf}{\zeta\left(\frac{1}{2}+it\right)}

\DeclareMathOperator{\im}{Im} 
\DeclareMathOperator{\cn}{cn}

\theoremstyle{definition}

\theoremstyle{remark}
\newtheorem{remark}[]{Remark}

\newtheorem*{mydef11}{{\bf Theorem 1}}

\newtheorem*{mydef12}{{\bf Theorem 2}}

\newtheorem*{mydef4}{{\bf Corollary}}

\newtheorem*{mydef91}{{\bf Formula1}}

\newtheorem*{mydef92}{{\bf Formula2}}

\newtheorem*{mydef93}{{\bf Formula3}}

\newtheorem*{mydef94}{{\bf Formula4}}

\numberwithin{equation}{section}

\begin{document}

\title{Jacob's ladders, class of crossbreeding conserving certain cell of meta-functional equations and existence of cancerous growth of that cell} 

\author{Jan Moser}

\address{Department of Mathematical Analysis and Numerical Mathematics, Comenius University, Mlynska Dolina M105, 842 48 Bratislava, SLOVAKIA}

\email{jan.mozer@fmph.uniba.sk}


\begin{abstract}
In this paper a set of 48 crossbreedings on certain cell of meta-functional equations preserving the cell is obtained. In opposite direction we presente an example of very complicated meta-functional equation not obtained in basic cell though it is generated by sequence of internal crossbreedings. Such a formula is called as \emph{cancerous growth} on the basic cell.  \\
\vspace{0.2cm} \ 

\noindent 
DEDICATED TO GREGOR MENDEL AND IT'S PEA-CROSSBREEDING

\end{abstract}
\maketitle

\section{Introduction} 

\subsection{}  

In the paper \cite{14}, eq. (1.10) we have considered the following set 
\be \label{1.1} 
\begin{split}
& \left\{\{\zeta(s),\Gamma(s),\cn(s,k), J_p(s)\}\right. \\ 
& \left. \{\Gamma(2s),\cn(2s,k),J_p(2s),\zeta(2s)\}\right. \\ 
& \left. \{\cn(3s,k),J_p(3s),\zeta(3s),\Gamma(3s)\}\right. \\ 
& \left. \{J_p(4s),\zeta(4s),\Gamma(4s),\cn(4s,k)\}\right. \\ 
& \left. \{\zeta(5s),\Gamma(5s),\cn(5s,k), J_p(5s)\}\right. \\ 
& \left. \{\Gamma(6s),\cn(6s,k),J_p(6s),\zeta(6s)\}\right. \\ 
& \left. \{\cn(7s,k),J_p(7s),\zeta(7s),\Gamma(7s)\}\right. \\ 
& \left. \{J_p(8s),\zeta(8s),\Gamma(8s),\cn(8s,k)\} \right. \\ 
& \left. \{\zeta(9s),\Gamma(9s),\cn(9s,k),J_p(9s)\} \right. \\ 
& \left. \hspace{2cm} \vdots  \right. \\ 
& \left. 
\right\}, \\ 
& s\in \mbb{C}\setminus\{N,P\},\ k^2\in (0,1) 
\end{split}
\ee 
of foursomes, where $\{N,P\}$ stands for the set of all zeros and all poles of the functions contained in (\ref{1.1}) for every admissible and fixed $k$. We have constructed, for example, the following set of exact meta-functional  equations on the subset of the first four sets in (\ref{1.1}) (for more information see \cite{14}, (3.8) -- (3.13)): 

\be \label{1.2} 
\begin{split}
& |\zeta(s_1^1)||\Gamma(s_2^1)||J_p(2s_3^2)|+|\zeta(2s_4^2)||\cn(s_3^1,k)|=\\ 
& = |\Gamma(2s_1^2)||\cn(s_3^1,k)||\cn(2s_2^2,k)|+|J_p(s_4^1)||J_p(2s_3^2)|, 
\end{split}
\ee 
\be \label{1.3} 
\begin{split}
	& |\zeta(s_1^1)||\zeta(3s_3^3)||\Gamma(s_2^1)|+|\Gamma(3s_4^3)||\cn(s_3^1,k)|=\\ 
	& = |\cn(s_3^1,k)||\cn(3s_1^3,k)||J_p(3s_2^3)|+|J_p(s_4^1)||\zeta(3s_3^3)|, 
\end{split}
\ee 
\be \label{1.4} 
\begin{split}
	& |\zeta(s_1^1)||\Gamma(s_2^1)||\Gamma(4s_3^4)|+|\cn(s_3^1,k)||\cn(4s_4^4,k)|=\\ 
	& = |\zeta(4s_2^4)||J_p(4s_1^4)||\cn(s_3^1,k)|+|J_p(s_4^1)||\Gamma(4s_3^4)|, 
\end{split}
\ee 
\be \label{1.5} 
\begin{split}
	& |\zeta(3s_3^3)||\Gamma(2s_1^2)||\cn(2s_2^2,k)|+|\Gamma(3s_4^3)||J_p(2s_3^2)|=\\ 
	& = |J_p(2s_3^2)||J_p(3s_2^3)||\cn(3s_1^3,k)|+|\zeta(2s_4^2)||\zeta(3s_3^3)|, 
\end{split}
\ee 
\be \label{1.6} 
\begin{split}
	& |\Gamma(2s_1^2)||\Gamma(4s_3^4)||\cn(2s_2^2,k)|+|J_p(2s_3^2)||\cn(4s_4^4,k)|=\\ 
	& = |J_p(2s_3^2)||J_p(4s_1^4)||\zeta(4s_2^4)|+|\zeta(2s_4^2)||\Gamma(4s_3^4)|, 
\end{split}
\ee 
\be \label{1.7} 
\begin{split}
	& |J_p(3s_2^3)||\Gamma(4s_3^4)||\cn(3s_1^3,k)|+|\zeta(3s_3^3)||\cn(4s_4^4,k)|=\\ 
	& = |\zeta(3s_3^3)||\zeta(4s_2^4)||J_p(4s_1^4)|+|\Gamma(3s_4^3)||\Gamma(4s_3^4)|, 
\end{split}
\ee 
i.e. there are the sets 
\bdis 
\Omega_l^m\subset\mbb{C},\ m,l=1,2,3,4 
\edis 
such that for their elements (variables) 
\bdis 
s_l^m\in\Omega_l^m 
\edis 
equations (\ref{1.2}) -- (\ref{1.7}) hold true. 

\begin{remark}
Let us remind that the equations in question are generated by the following mother formula (exact complete hybrid formula): 
\be \label{1.8} 
\tilde{Z}^2(\alpha_1^{1,1})\sin^2(\alpha_0^{1,1})+\tilde{Z}^2(\alpha_1^{2,1})\cos^2(\alpha_0^{2,1})=\tilde{Z}^2(\beta_1^1), 
\ee 
(see \cite{10}, (3.2), $k_1=k_2=1$ and $\frac{\pi}{4}\to \frac{\pi}{2}$ in our case), by this way: first this formula generates set of four transmutations (see \cite{14}, (3.3) -- (3.6)), and secondly the operation of crossbreeding applies on the last set and noticed equations follow. 
\end{remark} 

\subsection{} 

Let the symbol 
\be \label{1.9} 
[(1.2)\times(1.3)]|\zeta(s_1^1)|
\ee 
denote crossbreeding on the set $\{(1.2),(1.3)\}$ of the meta-functional equations (\ref{1.2}) and (\ref{1.3}) for the neutral factor $|\zeta(s_1^1)|$. Next, let $C_0$ denote the cell of exact meta-functional equations (\ref{1.2}) -- (\ref{1.7}) (elements of $C_0$). 

\begin{remark}
The cell $C_m$ is in this context generated by the transformation 
\be \label{1.10} 
g \to 4m+g,\ g=1,2,3,4,\ m\in\mbb{N}_0
\ee  
in $C_0$ as follows 
\be \label{1.11} 
\begin{split}
& |\zeta(s_1^1)|\to |\zeta[(4m+1)s_1^{4m+1}]|,\ |\Gamma(s_2^1)|\to |\Gamma[(4m+1)s_2^{4m+1}]|, \dots \\ 
& |\Gamma(4s_3^4)|\to |\Gamma[(4m+4)s_3^{4m+4}]|. 
\end{split}
\ee 
\end{remark} 

In this paper we show that there is a class of 48 crosbreedings of type (\ref{1.9}) conserving the cell $C_0$. For example, we have 
\be \label{1.12} 
[(1.2)\times(1.3)]|\zeta(s_1^1)|=(1.5)\in C_0. 
\ee 

\subsection{} 

Let us notice, on the contrary, there is a chain of such crossbreedings in $C_0$ that gives, for example, the following result: 
\be \label{1.13} 
\begin{split}
& |\Gamma(s_2^1)||\Gamma(3s_4^3)||\Gamma(4s_3^4)||\zeta(s_1^1)||\zeta(2s_4^2)||J_p(3s_2^3)||\cn(3s_1^3,k)|+ \\ 
& |\Gamma(s_2^1)|\underline{|\Gamma(3s_4^3)||\Gamma(3s_4^3)|}|\zeta(s_1^1)||\zeta(4s_2^4)||J_p(2s_3^2)||J_p(4s_1^4)|+ \\ 
& |\zeta(3s_3^3)||\zeta(4s_2^4)||\Gamma(2s_1^2)||\Gamma(3s_4^3)||J_p(s_4^1)||J_p(4s_1^4)||\cn(2s_2^2,k)|+ \\ 
& 
|\zeta(s_1^1)||\zeta(2s_4^2)||\zeta(3s_3^3)||\Gamma(s_2^1)||\Gamma(3s_4^3)||\cn(4s_4^4,k)|+ \\ 
& 
|\Gamma(2s_1^2)|\underline{|\Gamma(3s_4^3)||\Gamma(3s_4^3)|}|\cn(s_3^1,k)||\cn(2s_2^2,k)||\cn(4s_4^4,k)| = \\ 
& = 
|\Gamma(2s_1^2)||\Gamma(3s_4^3)||\Gamma(4s_3^4)||J_p(s_4^1)||J_p(3s_2^3)||\cn(2s_2^2,k)||\cn(3s_1^3,k)|+ \\ 
& 
|\Gamma(2s_1^2)|\underline{|\Gamma(3s_4^3)||\Gamma(3s_4^3)|}|\cn(s_3^1,k)||\cn(2s_2^2,k)||\zeta(4s_2^4)||J_p(4s_1^4)|+ \\ 
& |\zeta(s_1^1)||\zeta(2s_4^2)||\zeta(3s_3^3)||\zeta(4s_2^4)||\Gamma(s_2^1)||\Gamma(3s_4^3)||J_p(4s_1^4)|+ \\ 
& 
|\Gamma(2s_1^2)||\Gamma(3s_4^3)||\cn(2s_2^2,k)||\cn(4s_4^4)||\zeta(3s_3^3)||J_p(s_4^1)|+ \\ 
& 
|\Gamma(s_2^1)|\underline{|\Gamma(3s_4^3)||\Gamma(3s_4^3)|}|\zeta(s_1^1)||J_p(2s_3^2)||\cn(4s_4^4,k)|.
\end{split}
\ee 

Let the symbol $(3,2)\leftrightarrow(3,2)$ denote the type of elements of $C_0=\{(1.2)-(1.7)\}$. 

\begin{remark}
Then the meta-functional equation (\ref{1.13}) - a cancerous growth of the basic cell $C_0$ - is of type 
\bdis 
(7,7,7,6,6)\leftrightarrow(7,7,7,6,6). 
\edis  
The comparison of the norms of these types 
\bdis 
\begin{split}
& |(3,2)\leftrightarrow (3,2)|=10;\ |(\dots)|=2\times(3+2), \\ 
& |(7,7,7,6,6)\leftrightarrow(7,7,7,6,6)|=66 
\end{split}
\edis 
gives us good reason to name it as \emph{cancerous growth} of $C_0$. 
\end{remark} 

\begin{remark}
We notice explicitly that the meta-functional equation (\ref{1.13}) contains four squares $|\Gamma(3s_4^3)|^2$. 
\end{remark} 

\begin{remark}
This paper is again based on new notions and methods in the theory of the Riemann's zeta-function we have introduced in our series of 53 papers concerning Jacob's ladders. These can be found in arXiv [math.CA] starting with the paper \cite{1}. 
\end{remark} 

\section{On the structure of the mother formula} 

\subsection{} 

Let us remind that 
\be \label{2.1} 
\begin{split}
& \alpha_1^{1,1}=\vp_1^0(d_1),\ \alpha_0^{1,1}=\vp_1^1(d_1),\ d_1=d_1(U,\pi L;f_1), \\ 
& \alpha_1^{2,1}=\vp_1^0(d_2),\ \alpha_0^{2,1}=\vp_1^1(d_2),\ d_2=d_2(U,\pi L; f_2), \\ 
& \beta_1^1=\vp_1^0(e),\ e=e(U,\pi L), 
\end{split}
\ee 
(see \cite{6}, (4.5) -- (4.17)), and 
\be \label{2.2} 
\begin{split}
& f_1=f_1(t)=\sin^2t,\ f_2=f_2(t)=\cos^2t, \\ 
& t\in [\pi L,\pi L+U],\ 0<U<\frac{\pi}{2},\ L\in\mbb{N}. 
\end{split}
\ee 
Consequently, we have the following complicated composite functions 
\be \label{2.3} 
\begin{split}
& \tilde{Z}^2(\alpha_1^{1,1})=\tilde{Z}^2\{\vp_1^0[d_1(U,\pi L;f_1)]\}, \\ 
& \sin^2(\alpha_0^{1,1})=\sin^2\{\vp_1^1[d_1(U,\pi L;f_1)]\}, \\ 
& \tilde{Z}^2(\alpha_1^{2,1})=\tilde{Z}^2\{\vp_1^0[d_2(U,\pi L;f_2)]\}, \\ 
& \cos^2(\alpha_0^{2,1})=\cos^2\{\vp_1^1[d_2(U,\pi L;f_2)]\},\\ 
& \tilde{Z}^2(\beta_1^1)=\tilde{Z}^2\{\vp_1^0[e(U,\pi L)]\}, 
\end{split}
\ee 
where 
\be \label{2.4} 
\begin{split}
& \tilde{Z}^2(t)=\frac{{\rm d}\vp_1(t)}{{\rm d}t}=\frac{|\zf|^2}{\omega(t)}, \\ 
& \omega(t)=\left\{1+\mcal{O}\left(\frac{\ln\ln t}{\ln t}\right)\right\}\ln t, 
\end{split}
\ee
see \cite{4}, (6.1), (6.7), (7.7), (7.8) and (9.1). 

Now we have the following detailed form 
\be \label{2.5} 
\begin{split}
& \tilde{Z}^2\{\vp_1^0[d_1(U,\pi L;f_1)]\}\sin^2\{\vp_1^1[d_1(U,\pi L;f_1)]\}+ \\ 
& + \tilde{Z}^2\{\vp_1^0[d_2(U,\pi L;f_2)]\}\cos^2\{\vp_1^1[d_2(U,\pi L;f_2)]\}=\\
& = \tilde{Z}^2\{\vp_1^0[e(U,\pi L)]\} 
\end{split}
\ee 
of the exact complete hybrid formula that is a mother formula in this context for us. 

\begin{remark}
Namely, formula (\ref{2.5}) generates miscellaneous sets of its transmutations. Next, the operation of crossbreeding (see \cite{9} -- \cite{14}) applied on these sets generates new miscellaneous sets of exact meta-functional equations. 
\end{remark} 

\begin{remark}
Let us notice explicitly that the mother formula (\ref{2.5}) is based mainly on the new transcendental functions 
\bdis 
\vp_1(t),\ \frac{{\rm d}\vp_1(t)}{{\rm d}t},\ \mbox{i.e.}\ \Omega_l^m=\Omega_l^m[\vp_1,\frac{{\rm d}\vp_1}{{\rm d}t}], 
\edis 
(comp. subsection 1.1) that we have introduced in our papers. 
\end{remark} 

\section{Jacob's ladders}   

\subsection{} 

Let us remind that the Jacob's ladder 
\bdis 
\vp_1(t)=\frac 12\vp(t) 
\edis 
we have introduced in \cite{1} (see also \cite{4}), where the function $\vp(t)$ is arbitrary continuous solution of the nonlinear integral equation (also introduced in \cite{1}) 
\bdis 
\int_0^{\mu[x(T)]}Z^2(t)e^{-\frac{2}{x(T)}t}{\rm d}t=\int_0^T Z^2(t){\rm d}t 
\edis 
where 
\bdis 
\begin{split}
	& Z(t)=e^{i\vartheta(t)}\zf , \\ 
	& \vartheta(t)=-\frac t2\ln\pi+\im\left\{\ln\Gamma\left(\frac 14+i\frac t2\right)\right\}, 
\end{split}
\edis 
where each admissible function $\mu(y)$ generates the solution 
\bdis 
y=\vp(T;\mu)=\vp(T);\ \mu(y)\geq 7y\ln y. 
\edis 
We call the function $\vp_1(T)$ the Jacob's ladder as an analogue of the Jacob's dream in Chumash, Bereishis, 28:12. 

\subsection{} 

Next, let us remind that the classical Hardy-Littlewood integral (1918) 
\be \label{3.1} 
\int_0^T\left|\zf\right|^2{\rm d}t
\ee 
has the following expression 
\be \label{3.2} 
\int_0^T\left|\zf\right|^2{\rm d}t=T\ln T+(2c-1-\ln 2\pi)T+R(T), 
\ee  
with, for example, the Ingham's estimate of the error term 
\be \label{3.3} 
R(T)=\mcal{O}(\sqrt{T}\ln T),\ T\to\infty. 
\ee 
However, by the Good's $\Omega$-theorem (1977) we have that 
\be \label{3.4} 
R(T)=\Omega(T^{1/4}),\ T\to\infty. 
\ee 

\begin{remark}
It follows from (\ref{3.4}) that 
\be \label{3.5} 
\limsup_{T\to\infty}|R(T)|=+\infty, 
\ee 
that is every expression of the type (\ref{3.2}) of the Hardy-Littlewood integral possesses an unbounded error term at $T\to\infty$. 
\end{remark}

\subsection{} 

Under the circumstances (\ref{3.2}) and (\ref{3.5}) we have proved that the Hardy-Littlewood integral (\ref{3.1}) has an infinite set of other completely new almost exact representations expressed by the following 

\begin{mydef91}
\be \label{3.6} 
\begin{split}
& \int_0^T\left|\zf\right|^2{\rm d}t=\vp_1(T)\ln\{\vp_1(T)\}+\\ 
& + (c-\ln 2\pi)\vp_1(T)+c_0+\mcal{O}\left(\frac{\ln T}{T}\right),\ T\to\infty
\end{split}
\ee 
with the property (comp. (\ref{3.5})) 
\be \label{3.7} 
\lim_{T\to\infty}R_1(T)=\lim_{T\to\infty}\left\{\mcal{O}\left(\frac{\ln T}{T}\right)\right\}=0, 
\ee 
where $c$ is the Euler's constant and $c_0$ is the constant from the Titchmarsh-Kober-Atkinson formula. 
\end{mydef91}  

\begin{remark}
Comparison between (\ref{3.5}) and (\ref{3.7}) sufficiently characterizes level of exactness of our representation (\ref{3.6}). 
\end{remark} 

Now, let us remind other formulae demonstrating power of Jacob's ladder $\vp_1(t)$. 

\subsection{} 

First, we have obtained the following formula (see \cite{2}, (1.1)). 

\begin{mydef92}
\be \label{3.8} 
\begin{split}
& \int_T^{T+U}\left|\zeta\left(\frac 12+i\vp_1(t)\right)\right|^4\left|\zf\right|^2{\rm d}t\sim\frac{1}{2\pi^2}U\ln^5T,\\ 
& U=T^{7/8+2\epsilon},\ T\to\infty. 
\end{split}
\ee 
\end{mydef92}  

\begin{remark}
Our formula (\ref{3.8}) is the first asymptotic formula of the sixth order expression in $|\zeta|$ on the critical line $\sigma=\frac 12$ in the theory of the Riemann's zeta-function. 
\end{remark} 

\subsection{} 

Now, let 
\bdis 
S(t)=\frac 1\pi\arg\left\{\zf\right\},\ S_1(T)=\int_0^T S(t){\rm d}t, 
\edis  
where $\arg$ function is defined by the usual way. We have obtained the following formula (see \cite{3}, (5.4), (5.5)). 

\begin{mydef93}
\be \label{3.9} 
\begin{split} 
& \int_T^{T+U}	\left\{\arg\zeta\left(\frac 12+i\vp_1(t)\right)\right\}^{2k}\left|\zf\right|^2{\rm d}t=\frac{(2k)!}{k!2^{2k}}U\ln T(\ln\ln T)^k, 
\end{split} 
\ee 
where $T\to\infty$ and $k$ is fixed positive number. 
\end{mydef93}  

\begin{mydef94}
\be \label{3.10} 
\int_T^{T+U}\{S_1[\vp_1(t)]\}^{2k}\left|\zf\right|^2{\rm d}t=a_kU\ln T,\ T\to\infty. 
\ee 
\end{mydef94} 

\begin{remark}
Formulae (\ref{3.6}), (\ref{3.8}) -- (\ref{3.10}) are not accessible within the current methods in the theory of the Riemann's zeta-function. 
\end{remark} 

\section{Crossbreeding conserving a cell of meta-functional equations} 

\subsection{} 

Let the cell $C_0$ be the set (\ref{1.2}) -- (\ref{1.7}) of meta-functional equations (elements). 

We apply the crossbreeding of type $(A_1)$ on two different elements of $C_0$ in case each of those two elements contains once a choice neutral factor.  In this case the crossbreeding $(A_1)$ is simply elimination of the neutral factor. For example in the case 
\bdis 
[(1.2)\times (1.3)]|\zeta(s_1^1)|
\edis 
we have in the complete form 
\be \label{4.1} 
\begin{split}
& (1.2) \ \Rightarrow \\ 
& |\zeta(s_1^1)|=\frac{|\Gamma(2s_1^2)||\cn(s_3^1,k)||\cn(2s_2^2,k)|}{|\Gamma(s_2^1)||J_p(2s_3^2)|}+\frac{|J_p(s_4^1)||J_p(2s_3^2)|}{|J_p(2s_3^2)||\Gamma(s_2^1)|}- \\ 
& \frac{|\zeta(2s_4^2)||\cn(s_3^1,k)|}{|\Gamma(s_2^1)||J_p(2s_3^2)|}, \\ 
& (1.3) \ \Rightarrow \\ 
& |\zeta(s_1^1)|=\frac{|\cn(s_3^1,k)||\cn(3s_4^3,k)||J_p(3s_2^3)|}{|\zeta(3s_3^3)||\Gamma(s_2^1)|}+\frac{|J_p(s_4^1)||\zeta(3s_3^3)|}{|\zeta(3s_3^3)||\Gamma(s_2^1)|}- \\ 
& \frac{|\Gamma(3s_4^3)||\cn(s_3^1,k)|}{|\zeta(3s_3^3)||\Gamma(s_2^1)|}, 
\end{split}
\ee  
that gives the equation (\ref{1.5}). 

\begin{remark}
The fact is that all the corresponding moduli are positive. This property follows from our construction of exact $\zeta$-factorization formulae (see \cite{7}, \cite{8}). 
\end{remark} 

Next we use the abbreviation (comp. (\ref{1.9})) 

\bdis 
\begin{split}
& ([(1.2)\times(1.3)][|\zeta(s_1^1)|,|\Gamma(s_2^1)|,|J_p(s_4^1)|]=(1.5))\Leftrightarrow \\ 
& \{([(1.2)\times(1.3)]|\zeta(s_1^1)|=(1.5))\wedge ([(1.2)\times(1.3)]|\Gamma(s_2^1)|=(1.5))\wedge \\ 
& ([(1.2)\times(1.3)]|J_p(s_1^4)|=(1.5))\}. 
\end{split}
\edis 
The complete result is as follows. 

\begin{mydef11}
For the cell $C_0=\{(1.2) - (1.7)\}$ there is the following set of 36 crossbreeding of the type $(A_1)$: 
\be \label{4.2} 
\begin{split}
& [(1.2)\times(1.3)][|\zeta(s_1^1)|,|\Gamma(s_2^1)|,|J_p(s_4^1)|]=(1.5), \\
& [(1.2)\times(1.4)][|\zeta(s_1^1)|,|\Gamma(s_2^1)|,|J_p(s_4^1)|]=(1.6), \\ 
& [(1.2)\times(1.5)][|\zeta(2s_4^2)|,|\Gamma(2s_1^2)|,|\cn(2s_2^2,k)|]=(1.3), \\ 
& [(1.2)\times(1.6)][|\zeta(2s_4^2)|,|\Gamma(2s_1^2)|,|\cn(2s_2^2,k)|]=(1.4), \\ 
& [(1.3)\times(1.4)][|\zeta(s_1^1)|,|\Gamma(s_2^1)|,|J_p(s_4^1)|]=(1.7), \\  
& [(1.3)\times(1.5)][|\Gamma(3s_4^3)|,|\cn(3s_1^3,k)|,|J_p(3s_2^3)|]=(1.2), \\ 
& [(1.3)\times(1.7)][|\Gamma(3s_4^3)|,|\cn(3s_1^3)|,|J_p(3s_4^3)|]=(1.4), \\  
& [(1.4)\times(1.6)][|\cn(4s_4^4,k)|,|\zeta(4s_2^4)|,|J_p(4s_1^4)|]=(1.2), \\ 
& [(1.4)\times(1.7)][|\cn(4s_4^4,k)|,|\zeta(4s_2^4)|,|J_p(4s_1^4)|]=(1.3), \\ 
& [(1.5)\times(1.6)][|\Gamma(2s_1^2)|,|\cn(2s_2^2,k)|,|\zeta(2s_4^2)|]=(1.7), \\ 
& [(1.5)\times(1.7)][|\Gamma(3s_4^3)|,|\cn(3s_1^3,k)|,|J_p(3s_2^3)|]=(1.6), \\ 
& [(1.6)\times(1.7)][|\cn(4s_4^4,k)|,|J_p(4s_1^4)|,|\zeta(4s_2^4)|]=(1.5), 
\end{split}
\ee 
and 
\be \label{4.3} 
 [(1.2)\times(1.7)]=\emptyset, \ [(1.3)\times(1.6)]=\emptyset,\ [(1.4)\times(1.5)]=\emptyset, 
\ee 
(i.e. in these cases the corresponding two sets of moduli do not contain common modulus), and apart of these we have also 
\be \label{4.4} 
\begin{split}
& [(1.2)\times(1.3)][|\zeta(s_1^1)||\Gamma(s_2^1)|]=(1.5), \\ 
& [(1.2)\times(1.4)][|\zeta(s_1^1)||\Gamma(s_2^1)|]=(1.6), \\ 
& [(1.3)\times(1.4)][|\zeta(s_1^1)||\Gamma(s_2^1)|]=(1.7). 
\end{split}
\ee  
Here the neutral factor is the product of two moduli. 
\end{mydef11} 

\subsection{} 

We apply the crossbreeding of type $(A_2)$ on two different elements of $C_0$ in case each of those two elements contains twice a choice neutral factor. 

In this case, for example 
\bdis 
[(1.2)\times(1.3)]|\cn(s_3^1,k)|, 
\edis  
we have in the complete form 
\be \label{4.5} 
\begin{split}
& (1.2) \ \Rightarrow \\ 
& \frac{|\cn(s_3^1,k)|}{|J_p(2s_3^2)|}\{|\zeta(2s_4^2)|-|\Gamma(2s_1^2)||\cn(2s_2^2,k)|\}=
|J_p(s_4^1)|-|\zeta(s_1^1)||\Gamma(s_2^1)|, \\ 
& (1.3) \ \Rightarrow \\ 
& \frac{|\cn(s_3^1,k)|}{|\zeta(3s_3^3)|}\{|\Gamma(3s_4^3)|-|\cn(3s_1^3,k)||J_p(3s_2^3)|\}=
|J_p(s_4^1)|-|\zeta(s_1^1)||\Gamma(s_2^1)|,
\end{split}
\ee 
thus we obtain (\ref{1.5}) by equating left-hand sides in (\ref{4.5}). 

The complete result is as follows. 

\begin{mydef12} 
For the cell $C_0=\{(1.2) - (1.7)\}$ there is the following set of 12 crosbreedings of type $(A_2)$:
\begin{align} \label{4.6} 
& [(1.2)\times(1.3)]|\cn(s_3^1,k)|=(1.5), & [(1.3)\times(1.7)]|\zeta(3s_3^3)|=(1.4), \nonumber \\ 
& [(1.2)\times(1.4)]|\cn(s_3^1,k)|=(1.6), & [(1.4)\times(1.6)]|\Gamma(4s_3^4)|=(1.2), \nonumber \\ 
& [(1.2)\times(1.5)]|J_p(2s_3^2)|=(1.3), & [(1.4)\times(1.7)]|\Gamma(4s_3^4)|=(1.3), \nonumber \\ 
& [(1.2)\times(1.6)]|J_p(2s_3^2)|=(1.4), & [(1.5)\times(1.6)]|J_p(2s_3^2)|=(1.7), \\ 
& [(1.3)\times(1.4)]|\cn(s_3^1,k)|=(1.7), & [(1.5)\times(1.7)]|\zeta(3s_3^3)|=(1.6), \nonumber \\ 
& [(1.3)\times(1.5)]|\zeta(3s_3^3)|=(1.2), & [(1.6)\times(1.7)]|\Gamma(4s_3^4)|=(1.5). \nonumber 
\end{align}
\end{mydef12} 

Now, from (\ref{4.2}) and (\ref{4.6}) we have 

\begin{mydef4}
For the cell $C_0$ there is the set of 48 crossbreedings of types $(A_1)$ and $(A_2)$ respectively, and these are conserving that basic cell. In other words, we have the set of 48 transmutations of the mother formula (\ref{1.8}), (=(\ref{2.5})). 
\end{mydef4} 

\begin{remark}
If we make use transformations (\ref{1.10}), (\ref{1.11}) (see Remark 2) on the results (\ref{4.2}) -- (\ref{4.4}), (\ref{4.6}) in $C_0$ then we obtain corresponding results in $C_m,\ m\in\mbb{N}$. 
\end{remark}

I would like to thank Michal Demetrian for his moral support of my study of Jacob's ladders.

\end{document}